\documentclass[a4paper,draft,10pt]{amsart}
\usepackage{amssymb}
\usepackage{amscd} 

\theoremstyle{plain}
 
 \newtheorem{prop}{Proposition}[section]

\newtheorem{conj}{Conjecture}[section]
\theoremstyle{definition}
 
 \newtheorem{dfn}{Definition}[section]
\theoremstyle{remark}
 
 \numberwithin{equation}{section}

\renewcommand{\leq}{\leqslant}
\renewcommand{\geq}{\geqslant}

\title[Difference of Powers of Consecutive Primes Which Are Perfect Squares]{Difference of Powers of Consecutive Primes Which Are Perfect Squares}



\author[Ventullo]{\bfseries Alessandro Ventullo}

\email{alessandroventullo@gmail.com}

\begin{document}

\vspace{18mm} \setcounter{page}{1} \thispagestyle{empty}

\begin{abstract}
We investigate the consecutive primes $p$ and $q$ ($p>q$) for which there exists a pair of natural numbers $(x,y)$ such that $p^x-q^y$ is a perfect square and make some conjectures.
\end{abstract}

\maketitle

\section{Introduction}  

Mih\u{a}ilescu's theorem says that the only solution to the equation $$x^a-y^b=1$$ in natural numbers $x,y>0$ and $a,b>1$ is $(x,a,y,b)=(3,2,2,3)$. Drawing inspiration by this beautiful result, we want to try to modify the problem. Indeed, we are interested in $x,y$ as consecutive primes and we release the condition that $x^a-y^b$ is equal to $1$, but we require that this is a perfect square. In other terms, we are interested in the consecutive primes $p,q$ such that $$p^x-q^y=n^2,$$ where $x,y,n \in \mathbb{N}$. We begin our analysis tackling the problem from the most elementary cases. We begin with the following.

\begin{prop}
The only pairs of natural numbers $(x,y)$ such that $3^x-2^y$ is a perfect square are $(0,0),(1,1),(2,3),(3,1),(4,5)$.
\end{prop}
\begin{proof}
Let $n \in \mathbb{N}$ such that 
\begin{equation}\label{first-eq}
3^x-2^y=n^2.
\end{equation}
Clearly, if $x=0$, then $y=0$. Let $x>0$. We have three cases.
\begin{description}
\item[(i)] If $y \geq 2$, then $3^x \equiv n^2 \pmod{4}$, so $x$ must be even, i.e. $x=2k$ for some $k \in \mathbb{N}^*$. Therefore, equation \eqref{first-eq} becomes $$(3^k-n)(3^k+n)=2^y.$$ It follows that $3^k-n=2^a$ and $3^k+n=2^b$, where $a,b \in \mathbb{N}$ and $a+b=y$. Moreover, $a<b$ and adding these two equations, we get $$2\cdot3^k=2^a+2^b.$$ If $a=0$, the LHS is even and the RHS is odd, contradiction. If $a \geq 2$, we obtain $2\cdot3^k \equiv 0 \pmod{4}$, contradiction. So, $a=1$ and $b=y-1$ and we get $$3^k=1+2^{y-2}.$$ If $y=3$, we obtain $k=1$, i.e. $x=2$. If $y \geq 4$, then $3^k \equiv 1 \pmod{4}$, so $k$ is even and we can write $$(3^{\frac{k}{2}}-1)(3^{\frac{k}{2}}+1)=2^{y-2}.$$ Since $3^{\frac{k}{2}}+1$ and $3^{\frac{k}{2}}-1$ are powers of $2$ and their difference is $2$, we obtain $3^{\frac{k}{2}}+1=4$ and $3^{\frac{k}{2}}-1=2$, i.e. $k=2$, which gives $x=4$ and $y=5$. Therefore, we obtain the solutions $(x,y) \in \{(2,3),(4,5)\}$.
\item[(ii)] If $y=1$, then $3^x-2=n^2$. We have that $x$ must be odd, otherwise $n^2 \equiv -1 \pmod{4}$, contradiction. In the ring of integers $\mathbb{Z}[\sqrt{-2}]$, we have $$3^x=(n-\sqrt{-2})(n+\sqrt{-2}).$$
Let $d=(n-\sqrt{-2},n+\sqrt{-2})$. Clearly, $d \ | \ 2\sqrt{-2}$, so $N(d) \ | \ 8$. Since $3^x$ is odd, then its norm is odd and this implies that $N(d)=1$, i.e. $d=\pm 1$. So, $n-\sqrt{-2}$ and $n+\sqrt{-2}$ are coprime in $\mathbb{Z}[\sqrt{-2}]$. Since these two factors have the same norm and the only non trivial factorization (up to sign permutations) of $3$ with factors with the same norm is $3=(1-\sqrt{-2})(1+\sqrt{-2})$, then this forces 
$$\begin{array}{rcl} n-\sqrt{-2}&=&(1-\sqrt{-2})^x \\ n+\sqrt{-2}&=&(1+\sqrt{-2})^x  \end{array} \qquad \textrm{or} \qquad \begin{array}{rcl} n-\sqrt{-2}&=&(1+\sqrt{-2})^x \\ n+\sqrt{-2}&=&(1-\sqrt{-2})^x  \end{array}$$
Considering the first equation, we get $\textrm{Im}(n-\sqrt{-2})=\textrm{Im}((1-\sqrt{-2})^x)$, i.e. $$-\sqrt{2}=\left[-{x \choose 1}+2{x \choose 3}-4{x \choose 5}+\ldots+(-1)^{\frac{x+1}{2}}2^{\frac{x-1}{2}}{x \choose x}\right]\sqrt{2},$$ i.e. $$-1=-{x \choose 1}+2{x \choose 3}-4{x \choose 5}+\ldots+(-1)^{\frac{x+1}{2}}2^{\frac{x-1}{2}}{x \choose x}.$$ Let $$\displaystyle f(x)=\sum_{k \textrm{ odd}}^x {x \choose k} (-1)^{\frac{k+1}{2}}2^{\frac{k-1}{2}}$$ be defined on the odd natural numbers. An easy check shows that $f(x) \leq f(x+2)$ for any odd $x$. Since $f(5)=11$, then $x \in \{1,3\}$. An easy check shows that $f(1)=f(3)=-1$ and we get $(x,y) \in \{(1,1),(3,1)\}$. If we consider the second system of equations, we get $\textrm{Im}(n-\sqrt{-2})=\textrm{Im}((1+\sqrt{-2})^x)$ and proceeding as before, we obtain no solutions.
\item[(iii)] If $y=0$, then $3^x-1=n^2$, i.e. $3^x=n^2+1$. Since $x>0$ and $n^2 \equiv 0,1 \pmod{3}$, we get no solutions in this case. 
\end{description}
In conclusion, $(x,y) \in \{(0,0),(1,1),(2,3),(3,1),(4,5)\}$.
\end{proof}
Using the same ideas, we also obtain the following propositions.

\begin{prop}
The only pairs of natural numbers $(x,y)$ such that $5^x-3^y$ is a perfect square are $(0,0),(1,0)$.
\end{prop}
\begin{proof}
Let $n \in \mathbb{N}$ such that 
\begin{equation}\label{second-eq}
5^x-3^y=n^2.
\end{equation}
Clearly, if $x=0$, then $y=0$. Let $x>0$. We have three cases.
\begin{description}
\item[(i)] If $y \geq 2$, then $5^x \equiv n^2 \pmod{9}$, so $x=6k$ or $x=6k+4$, where $k \in \mathbb{N}^*$. In the first case, equation \eqref{second-eq} becomes $$(5^{3k}-n)(5^{3k}+n)=3^y.$$ It follows that $5^{3k}-n=3^a$ and $5^{3k}+n=3^b$, where $a,b \in \mathbb{N}$ and $a+b=y$. Moreover, $a<b$ and adding these two equations, we get $$2\cdot5^{3k}=3^a+3^b.$$ If $a>0$, the RHS is divisible by $3$, but the LHS is not divisible by $3$, so $a=0$ and $b=y$ and we get $$2\cdot5^{3k}=1+3^y.$$ But then $2\cdot5^{3k} \equiv 1 \pmod{9}$, contradiction. If $x=6k+4$, we get $$(5^{3k+2}-n)(5^{3k+2}+n)=3^y.$$ Reasoning as before we obtain a contradiction, so there are no solutions in this case.
\item[(ii)] If $y=1$, then $5^x-3=n^2$. Reducing modulo $4$, we obtain $n^2 \equiv 2 \pmod{4}$, so no solutions in this case.
\item[(iii)] If $y=0$, then $5^x-1=n^2$.  We have that $x$ must be odd, otherwise $(5^{\frac{x}{2}}-1)(5^{\frac{x}{2}}+1)=n^2$ and since the two factors on the LHS are coprime, they must be both perfect squares, contradiction. In the ring of integers $\mathbb{Z}[i]$, we have $$5^x=(n-i)(n+i).$$
Let $d=(n-i,n+i)$. Clearly, $d \ | \ 2i$, so $N(d) \ | \ 4$. Since $5^x$ is odd, then its norm is odd and this implies that $N(d)=1$, i.e. $d=\pm 1, \pm i$. So, $n-i$ and $n+i$ are coprime in $\mathbb{Z}[i]$. Since these two factors have the same norm and the only non trivial factorization (up to sign permutations) of $5$ with factors with the same norm is $5=(2-i)(2+i)$, then this forces 
$$\begin{array}{rcl} n-i&=&(2-i)^x \\ n+i&=&(2+i)^x  \end{array} \qquad \textrm{or} \qquad \begin{array}{rcl} n-i&=&(2+i)^x \\ n+i&=&(2-i)^x  \end{array}$$
Considering the first equation, we get $\textrm{Im}(n-i)=\textrm{Im}((2-i)^x)$, i.e. $$-1=-2^{x-1}{x \choose 1}+2^{x-3}{x \choose 3}-2^{x-5}{x \choose 5}+\ldots+(-1)^{\frac{x+1}{2}}{x \choose x}.$$ Let $$\displaystyle f(x)=\sum_{k \textrm{ odd}}^x {x \choose k} (-1)^{\frac{k+1}{2}} 2^{x-k}$$ be defined on the odd natural numbers. An easy check shows that $f(x) \geq f(x+2)$ for any odd $x$. Since $f(3)=-11$, then $x<3$, i.e. $x=1$. An easy check shows that $x=1$ works, so $(x,y)=(1,0)$. If we consider the second system of equations, we get $\textrm{Im}(n-i)=\textrm{Im}((2+i)^x)$ and proceeding as before, we obtain no solutions.
\end{description}
In conclusion, $(x,y) \in \{(0,0),(1,0)\}$.
\end{proof}

\begin{prop}
The only pair of natural numbers $(x,y)$ such that $7^x-5^y$ is a perfect square is $(0,0)$.
\end{prop}
\begin{proof}
Let $n \in \mathbb{N}$ such that 
\begin{equation}\label{third-eq}
7^x-5^y=n^2.
\end{equation}
Clearly, if $x=0$, then $y=0$. Let $x>0$. Reducing the equation modulo $4$, we obtain that $x$ must be even, i.e. $x=2k$ for some $k \in \mathbb{N}^*$. So, \eqref{third-eq} becomes $$(7^k-n)(7^k+n)=5^y.$$ It follows that $7^k-n=5^a$ and $7^k+n=5^b$, where $a,b \in \mathbb{N}$ and $a+b=y$. Moreover, $a<b$ and adding these two equations, we get $$2\cdot7^k=5^a+5^b.$$ If $a>0$, the RHS is divisible by $5$, but the LHS is not divisible by $5$, so $a=0$ and $b=y$ and we get $$2\cdot7^k=1+5^y.$$ If $y \geq 2$, then $2\cdot7^k \equiv 1 \pmod{25}$, contradiction. If $y=1$, then $2\cdot7^k=6$, contradiction. If $y=0$, then $7^k=1$, which has no solutions if $k \in \mathbb{N}^*$. So, $(x,y)=(0,0)$.
\end{proof}

\section{Generalizations and Conjectures}

From what we have seen, we can give generalizations and conjectures to the problem.

\begin{dfn}
We say that two primes $p$ and $q$ ($p>q$) are \emph{trivially squared} if $p^x-q^y$ is a perfect square implies that $(x,y)=(0,0)$. Otherwise, we say that two primes $p$ and $q$ ($p>q$) are \emph{nontrivially squared}.
\end{dfn}

In the last proposition we used only the congruence modulo $4$. So, if there are infinitely many twin primes $p,p-2$, where $p \equiv 3 \pmod{4}$, we obtain that there are infinitely many twin primes $p,p-2$ ($p>3$) which are trivially squared. We can go further and give the following conjecture.

\begin{conj}
There are infinitely many consecutive primes $p$ and $q$ ($p>q$) which are trivially squared.
\end{conj}

If we know that there are infinitely many primes $p$ of the form $n^2+1$, then if $(x,y)=(1,0)$ we obtain that $p^x-q^y$ is a perfect square for infinitely many consecutive primes $p$ and $q$.

\begin{conj}
There are infinitely many consecutive primes $p$ and $q$ ($p>q$) which are nontrivially squared.
\end{conj}

Observe that if the Landau's conjectures are true, then the two conjectures above would be true. Another observation of what we have done raises some other questions. For example, we can notice that the pairs of consecutive primes $(3,2)$, $(5,3)$ and $(7,5)$ yield a finite number of pairs of natural numbers $(x,y)$ such that $p^x-q^y$ is a perfect square. We state that this happens in general.

\begin{conj}
For any pair of consecutive primes $p$ and $q$ ($p>q$) there are only a finite number of pairs of natural numbers $(x,y)$ such that $p^x-q^y$ is a perfect square.
\end{conj}



\bibliographystyle{amsplain}

\end{document}